# Mathematical Incompleteness Results in First-Order Peano Arithmetic: A Revisionist View of the Early History

Saul A. Kripke

*Distinguished Professor in the Philosophy and Computer Science Programs at The Graduate Center, CUNY, New York City, USA*

The Saul Kripke Center, Room 7118, The Graduate Center, CUNY, 365 Fifth Avenue, New York, NY 10016, USA. kripkecenter@gc.cuny.edu. ORCiD: 0000-0001-7993-9456

# Mathematical Incompleteness Results in First-Order Peano Arithmetic: A Revisionist View of the Early History

*Abstract.* In the Handbook of Mathematical Logic, the Paris-Harrington variant of Ramsey's theorem is celebrated as the first result of a long 'search' for a purely mathematical incompleteness result in first-order arithmetic. This paper questions the existence of any such search and the status of the Paris-Harrington result as the first mathematical incompleteness result. In fact, I argue that Gentzen gave the first such result, and it that it was restated by Goodstein in a number-theoretic form.

## 1. The Received Account[1]

The received account of the subject can be given in two parts. One is the editor's note by Jon Barwise in the *Handbook of Mathematical Logic*, introducing the classic Paris-Harrington result:

---

[1] This paper represented my views about the history of the subject sometime after the *Handbook of Mathematical Logic*'s characterization of the significance of the Paris-Harrington result appeared, and shortly after the publication of Kirby and Paris (1982), when I looked up Goodstein's original paper. Nevertheless, I only drafted the present paper in 2007 (though I discussed it in a seminar given at Princeton University in 1998 and I remember discussing it with some logicians as well). I believe that what I call the "received account" was indeed the

Since 1931, the year Gödel's Incompleteness Theorems were published, mathematicians have been looking for a strictly mathematical example of an incompleteness in first-order Peano arithmetic, one which is mathematically simple and interesting and does not require the numerical coding of notions from logic. The first such examples were found in 1977 … (*Barwise 1977*, p. 1133)[2]

The allusion is of course to *Paris and Harrington 1977*, a simplification of earlier work of Paris.

Later, other such independence results were found, the best known being those of *Kirby and Paris 1982*. One result (their second involving 'hydras') is less widely quoted than the other, though the basic points of the hydra result are closely related to their more celebrated proof of the independence of a result of Goodstein. Here, to take a common statement of the result, I cite the exposition in Henle's textbook (*1986*, pp. 45-46). He calls a number, say 23, *written to the 'superbase' 2*, when all numbers greater than 2 are eliminated from its base-2 expression:

$$23 = 2^4 + 2^2 + 2^1 + 2^0 = 2^{2^2} + 2^2 + 2^1 + 2^0 \text{ (or } = 2^2 + 2^{2^2} + 2 + 1)$$

Similarly, for 514

$$514 = 2^9 + 2 = 2^{(2^{(2+1)}+1)} + 2$$

In other words, to write a number to the base *k* is to write it as a sum of powers of *k* (as a polynomial in *k*, with all coefficients <*k*). Suppose we apply the same

---

received account, in spite of the reactions of Kreisel and Gödel at the time. However, things seem to have changed for the better, as will be clear from what is said in the later part of the paper, even about Goodstein's work, let alone Gentzen's.

[2] The present writer is unaware of such a search. Was he isolated from the relevant mathematical community?

representation to the exponents (in the polynomial representation) and iterate the process. Then, it is said, Goodstein defined the sequences of a (now) well-known kind. Given any term of the sequence, the base involved is replaced by its successor and one is subtracted. In general, the next item in the sequence will evidently be enormously larger than its predecessor. Nevertheless, *Goodstein 1944* showed that the sequence eventually decreases, and indeed terminates at 0! (In other words, an infinite Goodstein sequence is impossible.)

> To quote Henle, Goodstein's theorem,
>
> … is remarkable in many ways. First, it is such a surprising statement that it is hard to believe it is true. Second, while the theorem is entirely about *finite* integers, Goodstein's proof uses *infinite* ordinals. Third, 37 years after Goodstein's proof appeared, L. Kirby and J. Paris proved that the use of infinite sets is actually *necessary*. That is, this is a theorem of arithmetic that can't be proved arithmetically, but *only* using the extra powers of set theory! (*1986*, p. 45)[3]

(I want to emphasize that I have cited Henle's textbook for expository purposes only, to give a reference in print. It expresses the common understanding in the logical community of the relevant history – I do not mean to criticize a particular author.)

---

[3] As in the previous footnote, Henle (*1986*, p. 48) describes a long search, really unknown to the present writer, for a proof of the statement in first-order PA, and claims that, finally, Kirby and Paris proved it impossible. As far as I am personally aware, the Goodstein paper was generally forgotten by the community of logicians until Kirby and Paris revived it. And, as we shall see below, most people have never looked up the original paper and misstate its result. Goodstein actually found a statement *equivalent* to $\varepsilon_0$-induction, and therefore independent.

## 2. The 'True' History

(a) First, it is hard for the present writer to see why it was not Gerhard Gentzen who found the first 'mathematical incompleteness' in Peano arithmetic in 1936. Only if 'logic' includes informal set theory (and remember that Cantor, who formulated the result in question, used no formal logic) would Gentzen's theorem – the independence of $\varepsilon_0$-induction with each ordinal $<\varepsilon_0$ represented in extended Cantor normal form (i.e., to the 'superbase' $\omega$) – fail to be a 'mathematical' independence result. Moreover, the statement is not one 'cooked up' in order to prove it independent, but is a natural existing result in the prior mathematical literature. (This is in contrast with the Paris-Harrington statement.)

The very same *Handbook* that pronounces that Paris and Harrington found the first mathematical incompleteness result in Peano arithmetic states in an earlier article by Helmut Schwichtenberg (see *Barwise 1977*, pp. 868 and 869) that $\varepsilon_0$-induction is an 'equivalent formulation of (*) having a clear mathematical meaning'. ((*) is the reflection principle for first-order arithmetic: $\exists x Der_S(x, \ulcorner \varphi \urcorner) \to \varphi$).

*Remark 1*. Somewhere in his voluminous writings, Kreisel also, writing about the Paris-Harrington theorem and the *Handbook*'s characterization, reacted that Gentzen's result is equally a 'mathematical' independence in first-order Peano Arithmetic. (Or so the present writer recalls.)[4]

---

[4] Also, of course, Kreisel (anticipated in unpublished work of Gödel) developed Gentzen's work into his 'no counterexample' interpretation of Peano Arithmetic, and also of fragments with restricted induction. For $\Pi_2^0$ statements the characterization can be given directly, and does not need a 'no counterexample' interpretation. His notion of ordinal recursion later developed into such things as the Hardy hierarchy and the Wainer hierarchy (extensions of the *Grzegorczyk* hierarchy).

*Remark 2.* In the oral tradition of the time, Gödel was said to have reacted, 'What could be more mathematical than a Diophantine equation?' Indeed, by the time of the *Handbook*, *Matiyasevich 1970*, building on the earlier work of *Davis, Putnam, and J. Robinson 1961*, had shown, in effect, that in any consistent recursively axiomatized system (in which some weak theory such as R is interpretable), there is a Diophantine equation that is undecidable. (These can be over the natural numbers – i.e. the non-negative integers – or over all the integers, if they are defined.)

There is little point in arguing that the result is somehow a disguised version of the original incompleteness theorem of Gödel, and hence involves the coding of metamathematical notions. The real distinction has to be between statements unprovable because of a result about *arbitrary* recursively axiomatized systems containing a bit of arithmetic, and one involving a specific weakness in first-order Peano arithmetic. The latter applies to the Gentzen statement and to the Paris-Harrington statement.

In this respect, Barwise's editorial remark is especially misleading. Gödel's original *1931* paper was not about first-order Peano Arithmetic, but about – as in its title – '*Principia Mathematica* and Related Systems'. (A system amounting to first-order Peano arithmetic is indeed mentioned, but not stressed.) In this respect, the result on Diophantine equations does the job quite well. But if one is looking at first-order PA specifically, the situation is otherwise.

(b) Perhaps the issue is a subjective one. For example, Fairtlough and Wainer in the *Handbook of Proof Theory*, following the usual view of the logical community, pronounced that Gentzen's results on $\varepsilon_0$-induction and corresponding results on fragments of first-order arithmetic are '"logic" independence results' (*Buss 1998*, p. 190), as opposed to the Paris-Harrington result, though at least they mention Gentzen in this connection so that the reader can form her own judgment.[5] I should add that in no way do I intend to denigrate the classic and beautiful Paris-Harrington result, which has influenced an enormous amount of fruitful work,[6] though in my opinion it should be called the first result in finite combinatorics, rather than the first 'mathematical' result, shown to be independent in first-order Peano Arithmetic. However, the truth of the Paris-Harrington statement is normally proved as a consequence of the infinite Ramsey theorem. It is very hard for the present writer to see why this result, and hence its consequences, is not just as much a set-theoretic result as $\varepsilon_0$-induction. In fact, I see even the finite results as set-theoretic, and were thus properly included in the *Handbook of Mathematical Logic* in that section, even though they would have been intelligible, and

---

[5] However, see the remarks below on the original Goodstein paper. I think no one denies that it is a number theoretic proposition, as Goodstein says. I have already stressed (see footnote 2 above) that in its original formulation it is *equivalent* to $\varepsilon_0$-induction, as Goodstein rightly stresses.

[6] Actually, I confess that the Paris-Harrington statement does seem to me have a certain defect of beauty. The original finite Ramsey theorem was proved for arbitrary finite sets. The cardinality of the set was all that mattered, not the identity of the elements. The Paris-Harrington statement required the identification of the elements with particular natural numbers.

However, in *Graham, Rothschild, and Spencer 1990* a theorem of Van der Waerden is called in the preface (p. 5) 'the central result of Ramsey's theory'. The theorem states that if the positive integers are partitioned into two classes, then at least one of the classes must contain arbitrarily long arithmetic progressions (see p. 29). It is obviously about the positive integers.

Certainly, and regardless of any such considerations, the statement in the text about how the Paris-Harrington result inspired a great deal of fruitful work still stands, my own work on 'fulfillability' included.

might have been discovered, before set theory became a common tool of mathematics. Of course, $\varepsilon_0$-induction is equivalent to the statement that there are no infinite descending chains (in extended Cantor normal form) of ordinals $<\varepsilon_0$. So stated, the result is a $\Pi_1^1$ statement and it is not statable in first-order PA. However, Gentzen showed that any proof of a contradiction in first-order PA would imply the existence of a primitive recursive infinite descending chain of such ordinals. Hence the statement that no such chains exist is $\Pi_2^0$. By Gödel's second incompleteness theorem, it is unprovable in PA.[7]

Exactly the same device is employed by *Kirby and Paris 1982* in their Theorem 2. They show that in their Hercules-hydra game every strategy is a winning strategy for Hercules (a $\Pi_1^1$ statement). Then they proved that the $\Pi_3^0$ statement that every recursive strategy is a winning strategy is unprovable in PA.[8] However, when Kirby and Paris say that they present 'perhaps the first [independence result in PA] which is, in an informal sense, purely number-theoretic in character (as opposed to metamathematical or combinatorial)' (*1982*, p. 285) what they say ignores Gentzen's statement (the basis of their own paper!), which is neither metamathematical nor combinatorial. However, the Gentzen statement is not, as *he himself formulated it*, number-theoretic in character.[9]

---

[7] Actually, Gentzen eventually showed that the independence of $\varepsilon_0$-induction could be shown directly, without recourse to Gödel's second incompleteness theorem (*Gentzen 1943*). The same is true of the Paris-Harrington theorem even though it is not presented that way in the *Handbook* proof.

[8] It seems to me that the result can be improved to get a $\Pi_2^0$ independence result, as in Gentzen and Goodstein. In fact, since Gentzen proved that an inconsistency in PA would give a primitive recursive infinite descending chain, so that the existence of such a chain is independent in PA, then its 'translation' into the hydra game ought to be independent in PA, not just in the fragment with $\Sigma_1^0$ induction, as Kirby and Paris seem to think.

Everyone agrees, however, that the Kirby-Paris result on the termination of Goodstein sequences is a mathematical, indeed purely number theoretic, result unprovable in PA. And it is in fact deduced from $\varepsilon_0$-induction. However, the community appears to have read Kirby and Paris a bit carelessly and largely ignored Goodstein. Already on the first page, Kirby and Paris rightly say, 'the first result of our paper is an improvement of a theorem of Goodstein' (*1982*, p. 285). Remembering that their first result is not the *truth* of the statement that the Goodstein sequences terminate, but rather its *independence* in PA, one should understand the right situation. Goodstein already proved a related independence result, weaker than that proved by Kirby and Paris (because it proves the independence of an even stronger statement). In fact, Kirby and Paris explicitly state two facts. First, that Gentzen 'showed that using transfinite induction on ordinals below $\varepsilon_0$ one can prove the consistency of P' (*1982*, p. 287). Second, that 'Goodstein proved the following: if $h: N \to N$ is a non-decreasing function, define an $h$-Goodstein sequence $b_0, b_1, \ldots$ by letting $b_{i+1}$ be the result of replacing every $h(i)$ in the base $h(i)$ representation of $b_i$ by $h(i+1)$, and subtracting 1. Then the statement "for every non-decreasing $h$, every $h$-Goodstein sequence eventually reaches 0" is equivalent to transfinite induction below $\varepsilon_0$' (*1982*, p. 287). So stated, the Goodstein proposition is even stronger than the statement usually quoted, and therefore even more surprising.

As Goodstein states it, the statement is $\Pi_1^1$ and therefore not statable in first-order PA. In fact, since Gentzen's argument actually showed that the existence of primitive recursive descending sequences of ordinals $< \varepsilon_0$ is independent of first-order PA, equivalently the existence of primitive recursive Goodstein sequences is independent of first-order PA. This last is a $\Pi_2^0$ statement, statable in first-order PA.

---

[9] Here I am referring to Goodstein's reformulation, as stressed below, which *is* purely number theoretic in character.

It is unclear to me whether Goodstein was aware of this latter fact. Probably he was not. All that I can find in his paper that is relevant is the vague assertion about decreasing sequences of ordinals (in extended Cantor normal form): 'For every constructively given sequence of ordinals the sequence $m_r$ is general recursive though not perhaps primitive recursive in every case' (*1944*, p. 34).

Various remarks ought to be made about Goodstein's paper. First, Goodstein himself expresses no recognition of how very surprising his result is. This is emphasized in the paper by Kirby and Paris, even for the weaker form that is the focus of their paper, where 'superbases' are raised by 1. (Henle, as quoted above, also emphasizes it.) Second, most articles on the subject, including the Wikipedia entry on 'Goodstein sequences' and the entry on the same subject in Wolfram Math World, state that Goodstein sequences are those represented by Kirby and Paris. In the Wikipedia entry on Goodstein sequences, Goodstein's *original* more general notion is called 'an extended Goodstein sequence'.

A further point is the following. Goodstein puts little emphasis on the fact that he has found a number-theoretic equivalent of Gentzen's incompleteness result. Rather, the emphasis is on the supposed 'finitist' character of his proof, at least of special cases, say up to $\omega^{\omega^\omega}$.

One must also mention a few other things. Gentzen's techniques can be developed not only to show the simple consistency but even the 1-consistency ($\Sigma_1$-correctness) of first-order PA. Also, the Gentzen proof (even of 1-consistency) can be extended to the case where a one place function letter is added. (Given paring functions, etc., this is equivalent to the system with infinitely many function letters of any number of variables.)

On this basis, Goodstein's formulation with arbitrary Goodstein sequences could be stated and proved independent. Moreover, on this basis, the Kreisel results characterizing the $\Pi_2^0$ statements provable in first-order PA (in terms of ordinal recursion $< \varepsilon_0$) and the $\Pi_1^1$ form (or 'no counterexample') interpretation of arbitrary statements provable in first-order PA can be proved (*Kreisel 1951* and *1952*). (The latter was already known to Gödel much earlier but was not published. See the introductory note to *Gödel 1938*, section 7, pp. 82-84.)

Finally, one should also mention that the paper by Ketonen and Solovay 1981 on the Paris-Harrington theorem bases everything on Gentzen's work on the unprovability of $\varepsilon_0$-induction, and the corollaries by Kreisel (and Gödel) on which functions are provably recursive in first-order Peano arithmetic. The same is also true not only of the paper by Kirby and Paris, but also of the shorter proof of their result by Cichon 1983.

## 3. Appendix

This was the situation in the early years. However, already in 2010, Craig Smorynski, reviewing a book by Peter Smith, stresses that the Paris-Harrington statement was *not* the first mathematical statement proved independent of first-order PA, contrary to the impression many people appear to have (*Smorynski 2010*, *Smith 2007*). Rather, he mentions Gentzen and the applications of his work by Kreisel. (Also, he observes that the Paris-Harrington result is an improvement of earlier work of Paris.) Even more important, Smith himself on his *Logic Matters* blog, remarking on Smorynski's review, says that several people have warned him similarly about Gentzen, etc., and that he is embarrassed that he did not take it out of later versions of his book. (Probably it has been taken out of the second edition.)

Most important, *Michael Rathjen 2015* published a highly precise and informative treatment of Goodstein's original paper, including an interesting recovery of Bernays's correspondence with Goodstein, due to Jan von Plato. Bernays was chosen as the referee of Goodstein's paper. He pointed out that the termination of arbitrary Goodstein sequences cannot be formulated in first-order PA and that one needs to add a function letter, and then the Gentzen consistency proof would still apply. Unfortunately, this led Goodstein to delete the explicit claim of an independence result and concentrate on the supposed 'finitist' character of his proof, at least up to certain ordinals $< \varepsilon_0$, for example $\omega^{\omega^\omega}$. Bernays had cautioned him about this also, but in this case he ignored Bernays's advice. He could have found an independence result in first-order PA itself using primitive recursive sequences; but as we saw above, he doesn't seem to be aware of this. Goodstein does vaguely talk of general recursive sequences. If he thought that this would be independent in first-order PA, he could have stated an independence result in first-order PA.

Rathjen gives a very careful treatment of the work of Goodstein, more strictly done than in the present paper or Goodstein's original. Rathjen, too, even in his abstract, shares my doubts as to whether there was really a long search 'for strictly mathematical examples of an incompleteness in PA' (see also footnote 1). He also says in the concluding sentence of his abstract, 'However, in relation to independence results, we think that both Gentzen and Goodstein are deserving of more credit'. This accords with the view that I have taken in the present paper. Rathjen also gives a treatment of what he calls 'special Goodstein sequences'. These are those where the 'superbase' increases by one, and it is these that are treated by Kirby and Paris.[10]

---

[10] In footnote 4, p. 231, *Rathjen 2015* says that 'around 1979, Diana Schmidt proved that Kruskal's theorem elementarily implies that the ordinal representation system for $\Gamma_0$ is well-founded [Schmidt! See *Naming and Necessity*, pp. 83-84]. She even wrote (p. 61) that she didn't know

As we have seen, even today (2020) the Wikipedia entry[11] on Goodstein's theorem states that Goodstein's sequences are only these and calls Goodstein's actual sequences 'extended Goodstein sequences'. (We have also seen that the Wolfram Math World entry also states that Goodstein sequences are those where superbase is increased by one.) However, the Wikipedia entry is aware of Gentzen's very early priority over Paris and Harrington for a mathematical independence proof, but does not appear to give Goodstein such a priority.[12,13]

**References**

Barwise, J. (ed.). 1977. *Handbook of Mathematical Logic*, Amsterdam: North-Holland.

Buss, S. (ed.). 1998. *Handbook of Proof Theory*, Amsterdam: Elsevier.

Cichon, E. A. 1983. 'A short proof of two recently discovered independence results using recursion theoretic methods', *Proceedings of the American Mathematical Society*, **87** (4), 704–6.

Davis, M., Putnam, H., and Robinson, J. 1961. 'The Decision Problem for Exponential Diophantine Equations', *Annals of Mathematics*, **74** (3), 425-36.

---

of any applications of her result to proof theory. This is quite surprising since in conjunction with proof-theoretic work of Feferman and Schütte from the 1960's it immediately implies the nowadays celebrated result that Kruskal's theorem is unprovable in predicative mathematics'.

[11] My editors have objected to my using quotes from Wikipedia. I hope I will be forgiven for doing so!

[12] Although it does say that Goodstein's statement is an early example of a statement demonstrably independent of first-order PA, it attributes everything to the paper by Kirby and Paris.

[13] I would like to thank Harold Teichmann and Yale Weiss for their editorial help. Special thanks go to Romina Padró for her invaluable help in producing the present and the original 2007 version of this paper. This paper has been completed with support from the Saul Kripke Center at the City University of New York, Graduate Center.